\documentclass[12pt]{article}
\usepackage{amssymb,graphics,amsmath}
\usepackage[all]{xy}
\newtheorem{definition}{Definition}[section]

\newtheorem{corollary}{Corollary}[section]

\newtheorem{theorem}{Theorem}[section]

\newtheorem{conjecture}{Conjecture}[section]

\def\Q{\mathbb{Q}}
\def\Z{\mathbb{Z}}

\def\what{\widehat}
\def\proof{{\bf Proof) }}
\begin{document}
\title{On the non-existence of L-space surgery structure}
\author{Motoo Tange\footnotemark}
\date{}
\footnotetext[1]{The author was supported by COE program of Mathematical Department of Osaka University.}
\maketitle 
\abstract
We exhibit homology spheres which never yield any lens spaces by integral
Dehn surgery by using Ozsv\'ath and Szab\'o's contact invariant.

\section{Introduction}
Let $Y$ be a closed oriented 3-manifold.
In this paper we write by $Y_r(K)$ the Dehn surgered manifold of a knot $K$ in $Y$ and with slope $r$.
Lens spaces can be obtained from the Dehn surgery of the unknot $U$
with slope $-p/q$, i.e. $L(p,q)=S^3_{-p/q}(U)$.

In general it is difficult to determine when a lens space can be obtained by an integral surgery of a non-trivial knot $K$ in $S^3$.
There are some well-known non-trivial knots in $S^3$ yielding lens spaces by integral surgeries, for example
torus knots, 2-cable knots of torus knots, and $(-2,3,7)$-pretzel knot and so on.

If we weaken the ambient space to any homology sphere,
we can construct many lens spaces by integral Dehn surgery.
For example in \cite{[1]} R. Fintushel and R. Stern have asserted that a lens space $L(p,q)$ is obtained by an 
integral Dehn surgery of a knot $K$ in a homology 3-sphere $Y$ if and only if $q=\pm x^2\bmod p$ holds for an integer $x$.
What kinds of homology spheres and the knots can realize for the solution $x$?
Are such homology spheres and knots restricted?
These questions are quite natural ones.

In the present paper we shall show that such homology spheres are restricted (see Theorem~\ref{introthe}).
We denote by $\overline{\Sigma(2,3,5)}$ the Poincar\'e homology sphere with the reversed-orientation as
the usual one.
\begin{theorem}
\label{introthe}
$\overline{\Sigma(2,3,5)}$ does not yield any lens spaces by any positive integral Dehn surgeries.
\end{theorem}
From contact topological view point we would like to consider this theorem.

We explain the original motivation for this paper.
In the author's previous paper \cite{[16]} many lens spaces have been constructed by positive 
Dehn surgeries on $\Sigma(2,3,5)$.
But from L-space homology spheres whose correction terms are neither $0$ nor $2$,
the author could not construct any lens spaces by any positive surgeries in some range of the order
$p$ of the 1st homology group and correction terms of homology spheres.
In particular from $\overline{\Sigma(2,3,5)}$ the author could not construct
any lens spaces in some range of $p$.
The notions of L-space and correction term shall be defined in Section~2.

On the other hand J.B. Etnyre and K.Honda in \cite{[15]} have shown that 
there do not exist any tight contact structures over $\overline{\Sigma(2,3,5)}$.
Hence I wondered if the two phenomena are related to each other.
This is main motivation of this paper.

In \cite{[16]} the author has studied examples yielding lens spaces by using the surgery formula 
between invariants of $Y$, $Y_0(K)$, and $Y_p(K)$.
Viewing the examples the author guesses if $Y\in \{\#^n\Sigma(2,3,5)\#^m\overline{\Sigma(2,3,5)}\}$ carries a positive lens surgery structure, then
$Y$ is $S^3$ or $\Sigma(2,3,5)$.

\section{Lens space (or L-space) surgery structure}
\label{structure}
We call a rational homology 3-sphere {\it L-space} if the Heegaard Floer homology of the manifold
is isomorphic to that of $S^3$.
The definition of Heegaard Floer homology is in \cite{[18]} and \cite{[19]}.
It is well-known that the set of L-spaces contains all spherical manifolds and some hyperbolic manifolds.
In order to generalize Theorem~\ref{introthe} to L-space surgery, we need
irreducibility of $Y-K$.
In fact the trivial $1$-surgery of $\overline{\Sigma(2,3,5)}$ gives rise to L-space clearly.

The coefficients of the homologies are always $\Z_2$ hence $HF^+(Y,\frak{s})$ is a $\Z_2[U]$-module,
where $U$ is the action to lower the degree of the homology by $2$.
When $Y$ is a rational homology sphere, the homology admits the
absolute $\Q$-grading as in \cite{[3]}.
The {\it correction term} $d(Y,\frak{s})$ is defined to be the minimal grading in the $HF^\infty(Y,\frak{s})$-part,
which is the image by the natural map $\pi_*:HF^\infty(Y,\frak{s})\to HF^+(Y,\frak{s})$
defined in \cite{[18]}.

\begin{definition}
Let $Y$ be a closed oriented 3-manifold.
We say that $Y$ carries positive (negative) L-surgery structure, if $p$ is a positive (negative)
integer and $Y_p(K)$ is an L-space for a null-homologous knot $K\subset Y$.
Moreover if $Y-K$ is irreducible, we say that $Y$ admits proper L-surgery structure.

In particular we say that $Y$ carries positive (negative) lens surgery structure if
$Y_p(K)$ is a lens space for a positive (negative) integer $p$.
\end{definition}
We suppose that any connected-sum component of $Y$ is not a lens space.
If $Y$ carries positive (or negative) lens surgery structure,
then obviously $Y$ carries positive (or negative) proper L-surgery structure.

$S^3$ carries both positive and negative lens surgery structure.
The homology sphere $\Sigma(2,3,5)$ carries positive lens surgery structure (see \cite{[16]}).
But for any homology sphere lens surgery structure does not always exist.
We shall prove the following in Section~\ref{mainpf}.
\begin{theorem}
\label{main}
Let $Y$ be an L-space homology sphere.
If $Y$ carries positive proper L-surgery structure, then $Y$ admits positive tight contact structure.
\end{theorem}
Let $Y$ be an L-space homology sphere.
If $Y_p(K)$ is an L-space, then $\widehat{HFK}(Y,K,g)\cong \Z_2$ holds by using the same method as \cite{[13]}.
It follows from this fact and Y. Ni's result in \cite{[14]} that if $Y-K$ is irreducible then $K$ is a fibered knot.
The tight contact structure above is the one induced from the open book decomposition along 
the knot $K$ by the method \cite{[17]} of W. Thurston and H. Winkelnkemper.

From Theorem~\ref{main} immediately we can prove the following non-existence result for $\overline{\Sigma(2,3,5)}$.
\begin{corollary}
\label{cor1}
$\overline{\Sigma(2,3,5)}$ does not carry positive proper L-surgery structure.
\end{corollary}
As we mentioned in Section~1, $\overline{\Sigma(2,3,5)}$ carries non-proper positive L-surgery structure.\\
\proof
Suppose that $\overline{\Sigma(2,3,5)}$ carries positive proper L-surgery structure.
From Theorem~\ref{main} $\overline{\Sigma(2,3,5)}$ must admit a positive tight contact structure.
But $\overline{\Sigma(2,3,5)}$ does not admit any positive tight contact structures by the work \cite{[15]} by J. B. Etnyre and K. Honda.
This is inconsistent.\\
\hfill$\Box$\\
{\bf Proof of Theorem~\ref{introthe})}
It follows from Corollary~\ref{cor1} and irreducibility of any lens space that this theorem is true.
\hfill$\Box$\\

\begin{corollary}
\label{cor2}
The homology sphere $\Sigma(2,3,5)\#\overline{\Sigma(2,3,5)}$ does not carry any proper L-surgery structure positive or negative.
\end{corollary}
\proof
Since $\Sigma(2,3,5)\#\overline{\Sigma(2,3,5)}$ does not have any tight contact structures positive or negative,
it follows from Theorem~\ref{main} that $\Sigma(2,3,5)\#\overline{\Sigma(2,3,5)}$ does not have
positive or negative L-surgery structures.\\
\hfill$\Box$\\
\section{Ozsv\'ath-Szab\'o's contact invariant}
\label{invariant}
Here we review the Ozsv\'ath-Szab\'o's contact invariant \cite{[20]}.
For a positive cooriented contact structure $\xi$ over a closed oriented
3-manifold $Y$
we will define the contact invariant $c(\xi)$ (originally defined in \cite{[20]}).
By E. Giroux's work \cite{[21]} each isotopy class of contact structures
exactly corresponds to an open book decomposition up to positive stabilization.
An open book decomposition is a triple $D=(Y,K,\pi)$, where
$K$ is a fibered knot in $Y$ and $\pi$ is the fibration map $Y-K\to S^1$.
The correspondence between an open book decomposition $D$ and a contact
structure is due to W. Thurston and H. Winkelnkemper \cite{[17]}.
We denote the correspondence as follows:
$$\{\text{open book decompositions}\}/\text{positive stabilization}\leftrightarrow
\{\text{contact structures}\}/\text{isotopy}$$
$$D\to \xi_D.$$

Let $(Y,\eta)$ be a contact structure and $(Y,K,\pi)$ the open book decomposition of $Y$ 
corresponding to $\eta$. 
We denote by the symbol $\frak{t}(\eta)$ the $\text{spin}^c$ structure associated with $\eta$.
For $Y_0(K)$ there is the canonical contact structure $\xi_0$ satisfying $c_1(\xi_0)[\hat{F}]=2g(\hat{F})-2$,
where $\hat{F}\subset Y_0(K)$ is the capped surface of a fiber $F$ of $\pi$.
Then the natural $\text{spin}^c$ cobordism $(-Y_0(K),\overline{\xi}_0)\to (-Y,\overline{\xi})$
induces $\widehat{HF}(-Y_0(K),\frak{t}(\overline{\xi}_0))\cong \Z_2\oplus \Z_2\to \widehat{HF}(-Y,\frak{t}(\overline{\xi}))$.
The invariant $c(\xi)$ is the image of the generator by the natural map
$HF^+(-Y_0(K),\frak{t}(\overline{\xi}_0))\cong \Z_2\to \widehat{HF}(-Y_0(K),\frak{t}(\overline{\xi}_0))$.
The main property used here is the following:
\begin{theorem}[\cite{[20]}]
\label{os}
If a positive contact structure $(Y,\xi)$ is overtwisted, then $c(\xi)=0$.
\end{theorem}
From this theorem $c(\xi)\neq 0$ implies tightness of $\xi$.
\section{Proof of Theorem~\ref{main}}
\label{mainpf}

\proof
Let $Y$ be an L-space homology sphere and $Y_p(K)$ an L-space.
By the result in \cite{[13]}, we have $\what{HFK}(Y,K,g)\cong \Z_2$.
Since $Y-K$ is irreducible, from the main theorem in \cite{[14]}, $K$ is fibered.
Hence we can build the contact structure over $Y$, which
is associated with the open book decomposition $D=(Y,K,\pi)$ induced from this fibration.
Consider the following surgery exact triangle.
\begin{equation}
\label{et}
\xymatrix{
HF^+(-Y)\ar[dr]^{F_3} &&\ar[ll]^-{F_1}  HF^+(-Y_0(K),Q^{-1}[i]) \\
  &  HF^+(-Y_p(K),[i])\ar[ur]^{F_2}& 
}
\end{equation}
Here $Q:\text{Spin}^c(-Y_0(K))\cong {\mathbb Z}\to\text{Spin}^c(-Y_p(K))$ is a map between two sets of
$\text{spin}^c$ structures
defined in \cite{[3]}.
Let $[i]$ denote the contact structure over $-Y_p(K)$ satisfying $Q(i)=[i]$.
If all $\text{spin}^c$-structures in $Q^{-1}[i]$ are $c_1\neq 0$, the map $F_2$ is the $0$-map.
Hence the map restricted of $F_1$ to the $\frak{t}(\overline{\xi}_0)$-component 
$$HF^+(-Y_0(K),\frak{t}(\overline{\xi}_0))=HF^+(-Y_0(K),1-g)\cong \Z_2 \to HF^+(-Y)$$
is injective.
Hence by definition of the contact invariant, $c({\xi_D})$ does not vanish.
From the Theorem~\ref{os} $\xi_D$ is a tight contact structure.

On the other hand in the case where $c_1(\frak{t}(\overline{\xi}_0))=0$ the genus of $K$ is one.
Then for non-zero $i$, $HF^+(-Y_0(K),i)$$=0$ and $HF_{\text{red}}(-Y_0(K),0)=0$.
The knot Floer homology of $K$ is $\widehat{HFK}(Y,K,i)\cong \Z_2$ for $i=0,\pm1$ and
$\cong 0$ for $i\neq 0,\pm1$.
In this case the tau invariant $\tau(-K)$ is $-1$ by the same method as \cite{[13]},
hence the contact invariant $c(\xi_D)$ does not vanish.
From Theorem~3.1 the contact structure $\xi_D$ is tight.
\hfill$\Box$\\

We call a knot $K$ in a homology sphere $Y$ {\it a lens space Berge knot} 
if an integral Dehn surgery of $K$ is a lens space and the dual knot $K'$ of $K$ is
the union of two arcs each of which is embedded in the meridian disk of
the genus one Heegaard decomposition (see Definition 1.7. in \cite{[13]}).

The author has also verified that many Brieskorn spheres appear 
as the homology spheres $Y$ yielding lens spaces.
Ozsv\'ath and Szab\'o have shown that any lens space Berge knot is fibered \cite{[13]}.
As a result many Brieskorn homology spheres carry proper L-surgery structure
with contact structures associated with the lens space Berge knots.
But I could not find proper L-surgery structures for Brieskorn homology spheres with the reversed orientation.
Here we state a conjecture which generalizes Corollary~\ref{cor1}.
\begin{conjecture}
Let $\Sigma(p,q,r)$ be any Brieskorn homology sphere.
Then $\Sigma(p,q,r)$ carries proper L-surgery structure but
$\overline{\Sigma(p,q,r)}$ does not carry proper L-surgery structure.
\end{conjecture}

\section*{Acknowledgments}
I thank Professor Ryosuke Yamamoto and T. Kadokami for their useful comments about contact structure
and Dehn surgery respectively.
It has inspired the author writing of the paper.
I also thank Professor M. Ue for his giving me many suggestions in Doctor course seminar.

 \noindent
 Motoo Tange\\
 Research Institute for Mathematical Sciences, \\
 Kyoto University, \\
 Sakyo-ku Kyoto-shi Kyoto-fu 606-8502, Japan. \\
 tange@kurims.kyoto-u.ac.jp

\end{document}